\documentclass[11pt]{amsart}
\usepackage{amsmath, amsthm, amsfonts, hyperref, graphicx, ifpdf, mathrsfs,color}
\usepackage{amssymb}
\usepackage{faktor}
\usepackage{bm}
\usepackage{enumerate}
\usepackage[utf8]{inputenc}
\usepackage{mathtools}
\usepackage{esint}
\usepackage{tikz-cd} 
\hypersetup{
	unicode=false,          
	pdftoolbar=true,        
	pdfmenubar=true,        
	pdffitwindow=false,     
	pdfstartview={FitH},    
	pdftitle={My title},    
	pdfauthor={Author},     
	pdfsubject={Subject},   
	pdfcreator={Creator},   
	pdfproducer={Producer}, 
	pdfkeywords={keywords}, 
	pdfnewwindow=true,      
	colorlinks=true,       
	linkcolor=blue,          
	citecolor=blue,        
	filecolor=magenta,      
	urlcolor=MidnightBlue          
}

\theoremstyle{plain}
\newcommand{\comment}[1]{{\color{red}{#1}}}

\newtheorem{theorem}{Theorem}[section]
\newtheorem{lem}[theorem]{Lemma}
\newtheorem{pro}[theorem]{Proposition}
\newtheorem{cor}[theorem]{Corollary}

\newtheorem{thx}{Theorem}

\newtheorem{theoremA}{Theorem}

\theoremstyle{definition}

\newtheorem{qua}[theorem]{Question}
\newtheorem{claim}[theorem]{Claim}

\theoremstyle{remark}
\newtheorem{rem}[theorem]{Remark}

\numberwithin{equation}{section}
\renewcommand{\qed}{\mbox{}\hfill\openbox}

\DeclareMathOperator{\PSL}{PSL(2,\C)}

\newcommand{\af}{almost Fuchsian}
\newcommand{\ap}{accidental parabolic}
\newcommand{\waf}{weakly almost Fuchsian}

\newcommand{\cs}{conformal structure}
\newcommand{\cvc}{convex core}
\newcommand{\fg}{fundamental group}

\newcommand{\gf}{geometrically finite}

\newcommand{\Hd}{Hausdorff dimension}
\newcommand{\hqd}{holomorphic quadratic differential}

\newcommand{\hym}{hyperbolic metric}
\newcommand{\htm}{hyperbolic three-manifold}

\newcommand{\kg}{Kleinian group}
\newcommand{\ksg}{Kleinian surface group}
\newcommand{\lmt}{limit set}

\newcommand{\maxp}{maximum principle}
\newcommand{\mi}{minimal immersion}

\newcommand{\ms}{minimal surface}

\newcommand{\pc}{principal curvature}
\newcommand{\qf}{quasi-Fuchsian}
\newcommand{\qfm}{quasi-Fuchsian manifold}
\newcommand{\RS}{Riemann surface}

\newcommand{\sff}{second fundamental form}

\newcommand{\TS}{Teichm\"{u}ller space}

\newcommand{\tm}{three-manifold}
\newcommand{\tg}{totally geodesic}

\newcommand{\woap}{without accidental parabolics}

\newcommand{\be}{\begin{equation}}
\newcommand{\ene}{\end{equation}}
\newcommand{\br}{\begin{rem}}
	\newcommand{\er}{\end{rem}}
\newcommand{\bl}{\begin{lem}}
	\newcommand{\bcor}{\begin{cor}}
		\newcommand{\ecor}{\end{cor}}
	\newcommand{\el}{\end{lem}}
\newcommand{\bd}{\begin{Def}}
	\newcommand{\ed}{\end{Def}}
\newcommand{\ben}{\begin{enumerate}}
	\newcommand{\een}{\end{enumerate}}
\newcommand{\bp}{\begin{proof}}
	\newcommand{\ep}{\end{proof}}
\newcommand{\bpo}{\begin{pro}}
	\newcommand{\epo}{\end{pro}}
\newcommand{\beq}{\begin{equation*}}
\newcommand{\eeq}{\end{equation*}}
\newcommand{\bear}{\begin{eqnarray*}}
	\newcommand{\eear}{\end{eqnarray*}}
\newcommand{\bt}{\begin{theorem}}
	\newcommand{\et}{\end{theorem}}

\newcommand{\HH}{\mathbb{H}^3}
\newcommand{\Mcal}{\mathcal{M}^3}

\newcommand{\Bcal}{\mathcal{B}}
\newcommand{\Gcal}{\mathcal{G}}

\newcommand{\Tcal}{\mathcal{T}_g(S)}
\newcommand{\R}{\mathbb{R}}
\newcommand{\C}{\mathbb{C}}
\newcommand{\AF}{\mathcal{AF}}
\newcommand{\QF}{\mathcal{QF}}

\renewcommand{\qed}{\mbox{}\hfill\openbox}
\numberwithin{equation}{section}

\allowdisplaybreaks

\def\XXint#1#2#3{{\setbox0=\hbox{$#1{#2#3}{\int}$}
		\vcenter{\hbox{$#2#3$}}\kern-.5\wd0}}

\makeatletter
\def\@citestyle{\m@th\upshape\mdseries}
\def\citeform#1{{\bfseries#1}}
\def\@cite#1#2{{%
		\@citestyle[\citeform{#1}\if@tempswa, #2\fi]}}
\@ifundefined{cite }{%
	\expandafter\let\csname cite \endcsname\cite
	\edef\cite{\@nx\protect\@xp\@nx\csname cite \endcsname}%
}{}
\makeatother

\begin{document}
	\parskip1ex
	\title[Minimal Surfaces in QF]
	{Beyond almost Fuchsian space}

	\author{Zheng Huang}
	\address[Z. ~H.]{Department of Mathematics, The City University of New York, Staten Island, NY 10314, USA}
	\address{The Graduate Center, The City University of New York, 365 Fifth Ave., New York, NY 10016, USA}
	\email{zheng.huang@csi.cuny.edu}
	
	\author{Ben Lowe}
	\address[B.~L.]{Department of Mathematics, University of Chicago, Chicago, IL. 60615, USA.}
	\email{loweb24@uchicago.edu}

	
	\begin{abstract}
		An {\af} manifold is a {\htm} of the type $S\times\R$ which admits a closed {\ms} (homotopic to $S \times \{0\}$) with the maximum {\pc} 
		$\lambda_0 <1$, while a {\waf} manifold is of the same type but it admits a closed {\ms} with $\lambda_0 \le1$. We first prove that any 
		{\waf} manifold is geometrically finite, and we construct a Canary-Storm type compactification for the {\waf} space. We use this to prove 
		uniform upper bounds on the volume of the {\cvc} and {\Hd} of the {\lmt}s of {\waf} manifolds, and to prove a gap theorem for the principal curvatures of minimal surfaces in hyperbolic 3-manifolds that fiber over the circle. We also give examples of {\qf} manifolds which admit unique closed {\ms}s without being {\waf}. 
	\end{abstract}
	
	\maketitle
	
	\section {Introduction}
	\subsection{Motivating Questions}
	Closed incompressible surfaces are fundamental in {\tm} theory. Thurston observed that a closed surface of {\pc}s less than 
	$1$ in magnitude is incompressible in a {\htm} and this was proved in \cite {Lei06}. In the setting of complete {\htm}s which are 
	diffeomorphic to $S \times \R$ ($S$ a closed surface of genus at least two), closed surfaces of small curvatures, especially when they 
	are also minimal, play an important role (see for instance \cite{Uhl83, Rub05, KS07, l21, CMN22}).

	There is a well-developed deformation theory for complete {\htm}s of the type $S \times \R$ (see for instance 
	\cite{Thu86, BB04, Min10, BCM12} and many others). For this class of {\htm}s {\woap}, non-degenerate ones are {\qf}. We denote the {\qf} space by $\QF$, 
	and the {\af} space, consisting of elements of {\qf} space that admit a closed {\it {\ms}} homotopic to $S\times \{0\}$ of {\pc}s less than one, by $\AF$. An {\af} manifold has many favorable 
	properties. For instance it admits an equidistant foliation by closed surfaces: if $M\in\AF$ and $\Sigma$ is the unique {\ms} in $M$, then 
	$M = \bigcup\limits_{r\in\R}\Sigma(r)$, where $\Sigma(r)$ is the surface with signed distance $r$ from $\Sigma$ (see calculations in 
	\cite{Uhl83, Eps84}). Many geometric quantities associated to $M$ are controlled in terms of $\lambda_0$, the maximum 
	of the {\pc} over $\Sigma$: minimal quasi-isometry constant for a map between $M$ and the corresponding Fuchsian manifold (\cite{Uhl83}), Teichm\"uller distance 
	between the conformal ends of $M$ (\cite {GHW10, Sep16}), volume of the {\cvc}, Hausdorff dimension of the limit set 
	(\cite{HW13}). Despite these estimates, some basic questions remain unanswered.  For instance, how bad can an {\af} manifold be, or can it degenerate? More specifically:
	\begin{qua}\label{qua1}

 Is there a bound on the volume of the convex core of an almost Fuchsian manifold $M \cong S \times \mathbb{R}$ depending only on the genus of $S$? 
	\end{qua}
	\noindent and similarly, 
	\begin{qua}\label{qua2}
 Is the Hausdorff dimension of the limit set of an almost Fuchsian manifold $M \cong S \times \mathbb{R}$ bounded above by $2-\epsilon$, for some $\epsilon>0$ depending only on the genus of $S$? 
	\end{qua}
	Furthermore, it is well-known that any {\af} manifold admits a unique closed {\ms}, while there are examples of {\qf} manifolds which admit multiple or 
	even arbitrarily many closed stable {\ms}s (\cite{And83, HW19}). One can thus ask whether having a unique closed {\ms} 
	characterizes the closure $\overline{\AF}$ of {\af} space within {\qf} space:
	\begin{qua}\label{qua3}
		Does there exist $M\in\QF\backslash\overline{\AF}$ such that $M$ admits a unique closed incompressible {\ms}?
	\end{qua}
	In this paper, we work to answer these questions and other related questions. Our approach is to construct a compactification of a wider class (we call 
	{\it {\waf} space}) of complete {\htm}s of the type $S \times \R$ where each element admits a closed {\ms} of the maximum {\pc} $\lambda_0 \le 1$. 
	This class was also considered by Uhlenbeck (\cite{Uhl83}).  
	\subsection{Notation and Terminology}
	We list a few notations we will frequently refer to in this paper. Throughout the paper $S$ is a closed surface of genus $g \ge 2$. 
	\ben[{\bf (1)}]
	\item
	The {\TS} $\Tcal$ is the space of {\cs}s on $S$, modulo biholomorphisms in the homotopy class of the identity.  Every {\cs} $\sigma \in \Tcal$ 
	on $S$ admits a unique {\hym} denoted by $g_\sigma$. 
	\item
	$\Mcal$ is the class of complete {\htm}s  diffeomorphic to $S\times\R$.
	\item
	For a closed incompressible surface $\Sigma\subset M\in \Mcal$, we always denote by $\lambda(\Sigma)$ its {\pc}s, and $\lambda_0$ the 
	maximum absolute value of a {\pc} over $\Sigma$.
	\item 
	$\Bcal\subset \Mcal$ is the subclass of $\Mcal$ such that each $B \in \Bcal$ admits a closed incompressible surface $\Sigma'$ 
	(homotopic to $S\times \{0\}$) with $|\lambda(\Sigma')|\le 1$, and there exists at least one point $p\in\Sigma'$ such that $|\lambda(p)|< 1$. 
	Note that we do not require $\Sigma'$ to be minimal here.
	\item
	$\QF\subset\Mcal$ is the space of {\qf} manifolds.
	\item
	$\AF\subset\QF$ is the space of {\af} manifolds. Each $M\in\AF$ admits a closed incompressible {\underline{\ms}} whose 
	{\pc}s are strictly less than one in magnitude, namely, $\lambda_0 < 1$. As a consequence of the {\maxp}, every {\af} manifold admits a unique 
	closed {\ms} (\cite{Uhl83}).
	\item
	$\Bcal_0$ is the subclass of $\Mcal$ such that each $B \in \Bcal_0$ admits a closed incompressible \underline{\ms} $\Sigma$ (diffeomorphic to $S$) 
	with $\lambda_0\le1$. We call such $B$ {\it \waf} and $\Bcal_0$ the {\it {\waf} space}. Note that by this definition, $\AF \subset \Bcal_0$.
	\een
	We know that $\Bcal_0\subset\Bcal$ because there must be some point $p\in \Sigma$ such that $\lambda(p) =0$. This is due to the fact that the 
	{\sff} of a {\mi} in a manifold of constant curvature is the real part of a {\hqd} (\cite{Hop89}), and any {\hqd} on a closed {\RS} of genus $g\ge 2$ 
	has exactly $4g-4$ zeros, counting multiplicity.
	\subsection{Main results}
	
	It has been an open question ([Theorem 3.3, \cite{Uhl83}]) whether $M$ must be  quasi-Fuchsian if $M\in \Mcal$ admits a closed incompressible {\ms} 
	$\Sigma$ such that $|\lambda(\Sigma)|\le 1$ and $\lambda_0 = 1$.  A partial answer was given in (\cite{San17}) where the author showed there are no doubly degenerate limits of {\af} groups. 
	Our first result is the following more general statement:
	
	\begin{thx}\label{main1}
		If $B$ is a complete {\htm} of the type $S\times\R$ and it admits a closed incompressible surface $\Sigma'$ 
		homotopic to $S$ (but not necessarily minimal) with $|\lambda(\Sigma')|\le 1$, if there exists at least one point $p\in\Sigma'$ such that $|\lambda(p)|< 1$, 
		and if the corresponding Kleinian group has no {\ap}s then $B$ is {\qf}.
	\end{thx}
	
	Note that while a {\qfm} does not have {\ap}s, apriori a {\waf} manifold may. See Section \ref{prelim} for definitions. As an immediate consequence, we have: 
	
	\begin{cor} \label{cor1}
		Any {\waf} manifold with no accidental parabolics is {\qf}. 
	\end{cor}
	
	When we do not assume away accidental parabolics, we are able to prove the following statement: 
	
	\begin{theoremA} \label{A2} 
		Every element of $\Bcal_0$ is geometrically finite, i.e., no element of  $\Bcal_0$ has a degenerate end.  
	\end{theoremA} 
	
	
	Rubinstein \cite{Rub05} gave an example that suggests that there might actually be elements of $\Bcal_0$ that have {\ap}s-- see Remark \ref{rubinsteinexample}.  It remains unclear whether $\Bcal_0 = \Bcal$. This is related to the following open question (\cite{Rub05}):  if $M \in \QF$ admits a closed incompressible surface
	(not necessarily minimal) of {\pc}s less than $1$ in magnitude, then is $M \in \AF$? 

	To derive rigidity properties for the {\waf} space, we construct a compactification $\overline{\Bcal_0}$ of $\Bcal_0$. This compactification is analogous to the 
	compactification of the space of {\ksg}s constructed by Canary-Storm (\cite{CS12}). The difference with their approach is that our compactification is defined in terms of data associated to the unique {\ms} in each element of $\Bcal_0$.  
	
	\begin{thx}\label{main2}
		There exists a compactification $\overline{\Bcal_0}$ of the (unmarked) {\waf} space $\Bcal_0$ that extends the Deligne-Mumford compactification 
		of moduli space of {\RS}s. 
	\end{thx}
	The points at infinity that we add correspond to disjoint unions of cusped {\waf} manifolds. 
	By working with {\ms}s our approach is adapted to the applications below, and seems to require less in the way of structural results about {\ksg}s, although we work in a more 
	specialized setting. 
	
	Theorem \ref{main2} allows us to answer Question ~\ref{qua1} and Question ~\ref{qua2}. Let $S$ be a closed surface.  
	\begin{cor}\label{cor3}
		There are $\epsilon>0$ and $L>0$ depending only on the genus of $S$ such that for every $M$ diffeomorphic to $S \times \mathbb{R}$ contained in $\Bcal_0$, the {\Hd} of the {\lmt} of $M$ is at most 
		$2-\epsilon$, and the volume of the {\cvc} of $M$ is at most $L$.  	
	\end{cor}

	Another application is to prove a gap theorem for the maximum principal curvature of minimal surfaces in doubly degenerate {\htm}. These arise, for instance, as infinite cyclic covering spaces of closed hyperbolic 3-manifolds that fiber over the circle. 
	  
	\begin{thx} \label{main4} 
For every $g \ge 2$ and $\rho_0 > 0$, there is an $\epsilon$ depending on $g$ and $\rho_0$ such that if a doubly degenerate hyperbolic 3-manifold homeomorphic to the product of a closed surface of genus $g$ with $\mathbb{R}$ has injectivity radius at least $\rho_0$, then any embedded {\ms} isotopic to the fiber has maximum principal curvature larger than $1+\epsilon$.
		\end{thx} 
	
This recovers a result of Breslin \cite{Bre11}, and will follow from Theorem \ref{alggeo} in \S 5.  Informally speaking, Theorem \ref{alggeo} states that in order for a sequence of doubly degenerate $M_n$ to contain stable minimal surfaces $\Sigma_n$  with maximum principal curvatures tending to 1,  algebraic and geometric limits must fail to agree along every subsequence.  McMullen (\cite{McM99}) gives a general criterion for algebraic and geometric limits to agree. His condition is that, for a sequence of hyperbolic elements approaching an accidental parabolic element in the algebraic limit,  the square of the imaginary part of the complex length divided by the real part of the complex length tends to zero. One could thus also formulate a minimal surface principal curvature gap theorem for Kleinian surface groups satisfying McMullen's criterion that generalizes Theorem \ref{main4}.  It would be interesting to determine the extent to which Theorem \ref{main4} is true without an assumption on the injectivity radius.  

\begin{rem} 
	Farre-Vargas-Pallete (\cite{FVP22}) recently obtained results similar to Theorem \ref{main4} for certain sequences of hyperbolic mapping tori. Their results are complementary to ours, as they prove a principal curvature gap theorem in the case where algebraic and geometric limits differ. Their approach is based on an analysis of how minimal surfaces interact with curves of short length in the ambient hyperbolic 3-manifold. A key step is to use the presence of cusps in the limit to produce embedded horocycles on the minimal surfaces that they consider (compare Remark \ref{rubinsteinexample}.)  
	 
\end{rem} 
	
	We also answer Question ~\ref{qua3} by constructing examples of {\qfm}s which are not {\waf} but that each admit a unique closed  {\ms}. 
	\begin{thx}\label{main3}
There exist $M \in \QF\backslash\Bcal_0$ that contain a unique closed {\ms} $\Sigma$ such that its maximum {\pc} $\lambda_0(\Sigma)$ is greater than one.  
	\end{thx}

	\subsection{Outline of the paper} After reviewing some relevant preliminary facts about {\ksg}s and {\af} manifolds in \S 2, we will prove our main results in 
	Sections \S 3, \S 4 and \S 5. In particular, Theorems A and A1 are proved in \S 3 and the compactification of the {\waf} space is constructed in \S 4, which proves 
	Theorem B and deduces some applications of the compactness result. Finally in \S 5 we prove Theorems C and D.
	
	\subsection{Acknowledgements} 
	We are grateful to Richard Canary, Christopher Leininger and Andrew Yarmola for helpful correspondence. The second author is especially grateful to Richard Canary for his generosity in responding to many emails.  We thank James Farre and Franco Vargas-Pallete for finding a mistake in the proof of a previous version of Theorem \ref{main4}.  Finally we also wish to thank an anonymous referee for his/her willingness to read the paper carefully and for many insightful comments. The research of Z. H. was supported by a PSC-CUNY grant and the research of B. L. was supported by an NSF grant DMS-2202830.
	
	\section{Preliminaries} \label{prelim} 
	\subsection{Kleinian Surface Groups} \label{prelimksg}
	A {\kg} $\Gamma$ is a discrete subgroup of $\PSL$, the orientation preserving isometry group of $\HH$. Any complete {\htm} can be written as	$\HH/\Gamma$ for some Kleinian group $\Gamma$. We will work with the Poincar\'e ball model for $\HH$ and denote by $S^2_\infty$ the sphere at infinity.  The \textit{limit set}
	$\Lambda_{\Gamma} \subset S^2_\infty$ of $\Gamma$ is the set of accumulation points on $S^2_\infty$ of the orbits $\Gamma x$ of $x \in \HH$. 
	The \textit{domain of discontinuity} of $\Gamma$ is defined to be $\Omega_{\Gamma} = S^2_\infty\backslash\Lambda_{\Gamma}$. In this paper we restrict ourselves to 
	the case that $\Gamma$ is a {\ksg}, namely, $\Gamma$ is isomorphic to the {\fg} of a hyperbolic surface $S$ of finite area, which we also assume to be orientable. Equivalently $\Gamma$ can be viewed as the image of a discrete and faithful representation from $\pi_1(S)$ to $\PSL$ . In the case that $S$ is a closed surface and $\Gamma$ contains no parabolics, it is a result of Thurston  and Bonahon (\cite{Bon86}) that the quotient $\HH/ \Gamma$ is diffeomorphic to $S \times\R$.

	The convex hull of $\Gamma$ is the smallest non-empty closed convex subset of $\HH$ invariant under $\Gamma$, and its quotient, $C(M)$, by $\Gamma$ is the 
	{\it {\cvc}} of the {\htm} $M = \HH/\Gamma$. We call $\Gamma$ {\it geometrically finite} if the volume of $C(M)$ is finite. Otherwise we call it 
	{\it geometrically infinite}.

	When the {\lmt} of $\Gamma$ is a round circle, it spans a {\tg} hyperbolic disk inside $\HH$, and $\Gamma$ is called a Fuchsian group and $M = \HH/\Gamma$ a 
	Fuchsian manifold. In fact it is a warped product of a hyperbolic surface with the real line, of the following explicit metric expression: $dt^2 + \cosh^2(t)g_\sigma$, 
	where $g_\sigma$ is the {\hym} on $S$. When $\Lambda_{\Gamma}$ is homeomorphic to a circle, we call $\Gamma$ (and the quotient manifold) {\it \qf}. Analytically a {\qfm} 
	is a quasi-conformal deformation of a Fuchsian manifold, and hence {\gf}, and the {\cvc} $C(M)$ is homeomorphic to $S \times [0,1]$.

	An element of $\pi_1(S)$ is non-peripheral if no loop representing it can be homotoped into a puncture.  A parabolic element $\gamma \in \Gamma$ is an {\it {\ap}} if it is the image of a non-peripheral loop in $\pi_1(S)$ under a representation of which $\Gamma$ is the image.    When a {\ksg} $\Gamma$ does not contain any {\ap}s, then (\cite{Thu79,Bon86}) 
	$\HH /  \Gamma$ is either {\qf}, or $C(M)$ is homeomorphic to $S\times [0,\infty)$ ({\it simply degenerate}) or $C(M)$ is homeomorphic to 
	$S\times (-\infty, \infty)$ ({\it doubly degenerate}).  
	
	A sequence of Kleinian groups $\Gamma_i$ converges to $\Gamma$ \textit{algebraically} if there are isomorphisms $\Gamma \rightarrow \Gamma_i$ that converge to the restriction to $\Gamma$ of the identity map $\PSL \rightarrow \PSL$.  In practice the following equivalent definition \cite{McM96}[3.1] will be easier to use: $\Gamma_i$ converges algebraically to $\Gamma$ if there are smooth homotopy equivalences $\mathbb{H}^3 / \Gamma \rightarrow \mathbb{H}^3/ \Gamma_i$, compatible with the respective markings by $\Gamma$ and $\Gamma_i$, whose restrictions to any fixed compact set converge $C^\infty$ to local isometries. 
 
 The sequence $\Gamma_i$ converges \textit{geometrically} to $\Gamma$ if there are choices of baseframes $\omega_i$ for $\mathbb{H}^3 \backslash \Gamma_i$ and $\omega$ for $\mathbb{H}^3 \backslash \Gamma$ and a sequence of balls $B_i \subset \mathbb{H}^3$ that exhaust $\mathbb{H}^3$ such that the center of $B_i$ projects to the basepoint of $\omega_i$ and  each $B_i/\Gamma_i$ can be mapped $k_i$-quasi-isometrically onto a subspace of $\mathbb{H}^3/\Gamma$ by differentiable maps $F_i$ with  $k_i \rightarrow 1$ as $i\rightarrow \infty$ and $DF_i(\omega_i)=\omega$ for all $i$.   We say $\Gamma_i$ converges \textit{strongly} to $\Gamma$ if it converges both algebraically and geometrically to $\Gamma$.

Finally we will need the following folklore lemma, which is well-known to the experts. Related statements were used in the chapter 9 of Thurston notes \cite{thurston2022geometry}. We provide a short proof here. 

\begin{lem} \label{density}
Let $M\cong S \times \mathbb{R}$ be a quasi-Fuchsian manifold for $S$ a closed surface of genus $g\geq 2$. Then there exists a universal constant $C_0$ so that every point in the convex core of $M$ is at a distance of at most $C_0$ from some closed geodesic in $M$.  
\end{lem}

\begin{proof}

The quasi-Fuchsian manifold $M$ is quasi-isometric to a Fuchsian manifold $F\cong S \times \mathbb{R}$, which  implies that there is a $\pi_1(S)$-equivariant homeomorphism between the round circle limit set $C$ of $F$ and the quasicircle limit set $L$ of $M$ \cite{McM96}[Theorem 2.5].  Note that pairs of points in $C$ joined by geodesics in $\mathbb{H}^3$ that project to closed geodesics in $F$ are dense in $C\times C$.  This is because $F$ is Fuchsian, and closed geodesics are dense in the unit tangent bundle of any closed hyperbolic surface.  Therefore pairs of points in $L$ joined by geodesics in $\mathbb{H}^3$ that project to closed geodesics in $M$ are dense in $L\times L$. This implies that any geodesic in $\mathbb{H}^3$ with endpoints in $L$ can be approximated as $C^1$-close as desired on compact sets of $\mathbb{H}^3$ by lifts of closed geodesics in $M$.  

The convex core of $M$ is the union of all ideal tetrahedra in $\mathbb{H}^3$ with vertices in the limit set of $M$ (this can be seen using the Klein ball model for $\mathbb{H}^3$, in which straight lines are  geodesics in the hyperbolic metric, and the fact that the analogous statement is true in Euclidean space, see \cite{rivin09}.)   The conclusion of the lemma then follows from the fact that any ideal tetrahedron in $\mathbb{H}^3$ is contained in the $C_0$-neighborhood of its edges, for some universal constant $C_0$.  

 
\end{proof}

	
	\subsection{Almost Fuchsian Manifolds} A {\qfm} is topologically $S\times \R$, but it can be very complicated geometrically. 
	Uhlenbeck and others introduced techniques of {\ms}s of small {\pc}s to study {\ksg}s. In particular Uhlenbeck (\cite{Uhl83}) defined the following subclass: a {\qfm} 
	$M$ is {\it \af} if it admits a closed {\it \ms} $S$ whose maximum {\pc} satisfies $\lambda_0 < 1$. A natural question is how far it is from the almost Fuchsian space to the boundary of the {\qf} space in the deformation space of {\ksg}s. Uhlenbeck was able to derive an explicit formula for the {\hym} of an {\af} manifold $M$ in terms of the {\cs} of the 
	unique {\ms} $S$ in $M$ and its {\sff}.  She accomplished this by showing that the normal exponential map at the unique minimal surface gave global coordinates 
	on $M$.  This enabled her to prove that the quasi-isometry constant for a quasi-isometry between an {\af} manifold and a Fuchsian manifold is bounded 
	above by $\frac{1+\lambda_0}{1-\lambda_0}$. Note that if $\lambda_0$  is allowed to be $1$, this estimate yields no information.

	Among the invariants of a {\qfm}, the volume of the {\cvc} and the {\Hd} of $\Lambda_{\Gamma}$ are particularly important. It is well-known that the 
	{\Hd} of the limit set defines an analytic function on {\qf} space (\cite{Rue82}), and it is valued between $1$ (which it attains exactly when $\Gamma$ is Fuchsian) and $2$ (\cite{Bow79}). 
	When $\Gamma$ is in the class of {\af}, we know the {\Hd} is bounded from above by $1+\lambda_0^2$, and the volume of the {\cvc} is bounded from above by 
	$4\pi(g-1)(\frac{\lambda_0}{1-\lambda_0^2}+ \frac12\ln\frac{1+\lambda_0}{1-\lambda_0})$ (\cite{HW13}). Neither estimate gives any information if $\lambda_0$ 
	is allowed to tend to $1$.

	Both \cite{Uhl83} and \cite{Eps86}  considered the case of {\it \waf}, namely, allowing $\lambda_0 = 1$. There exist global coordinates for a {\waf} manifold via the hyperbolic Gauss map from the minimal surface with principal curvatures smaller than or equal to 1 (see \cite{Uhl83}[pg. 161].) This is the key 
	geometric property of {\waf} manifolds that we will utilize.

\section{Weakly Almost Fuchsian is Geometrically Finite} \label{firststep} 
In this section we will prove Theorem A and  Theorem A1. Theorem A1 is the first and most important step in constructing the compactification in Theorem B.  Our arguments build on ideas of Sanders \cite{San17}.  

 \noindent \textit{Proof of Theorem A:}

Let $M \in \Bcal$, then $M\cong S \times \mathbb{R}$ contains a surface $\Sigma$ homotopic to $S \times \{0\}$ such that $|\lambda(\Sigma)|\le 1$ and this inequality is strict at at least one point. We will show that neither of the ends of $M$ is simply degenerate.  First we choose a properly embedded disk $\tilde{\Sigma}$ lifting $\Sigma$ to the universal cover 
$\tilde{M} \cong \HH$. The normal exponential map 
\beq
\eta: \tilde{\Sigma} \times \mathbb{R} \rightarrow \HH  
\eeq
is a diffeomorphism by, for example, the arguments in \cite{Eps86}. The forward (+) and backwards (-) Gauss maps 
\beq
\Gcal_{\tilde{\Sigma}}^{\pm}: \tilde{\Sigma} \rightarrow \partial_\infty \HH
\eeq
associated to $\tilde{\Sigma}$  are defined by 
\[
\Gcal_{\tilde{\Sigma}}^{\pm}(p) = \lim_{t \to \pm \infty} \eta(p,t), 
\]     
\noindent where we understand these limits to lie in the sphere at infinity in the Poincar\'e ball model for $\HH$.

Let $U$ be a small disk in $\tilde{\Sigma}$ such that $|\lambda(\tilde{\Sigma})| < 1$ on the closure of $U$, and so that $U$ is disjoint from all non-identity translates of itself under the covering action of $\pi_1(\Sigma)$.  It follows from Epstein \cite{Eps86}[Equation (5.5)] 
that $\Gcal_{\tilde{\Sigma}}^{\pm}|_{U}$ has quasi-conformal dilatation bounded from above by $\frac{2}{\epsilon}$, where 
\be\label{epsi}
\epsilon = \min_{i=1,2} \{1-|\kappa_i|\}, 
\ene
and $\kappa_1$ and $\kappa_2$ are the supremal and infimal {\pc}s of $\tilde{\Sigma}$ on $U$. We take some point on $U$ to be the origin in the Poincare ball model of $\mathbb{H}^3$, and give $\partial_{\infty} \mathbb{H}^3$ the corresponding unit sphere metric. The generalization of the Koebe $\frac14$-Theorem 
proved by Astala and Gehring (\cite{AG85}) implies that the image of $U$ under $\Gcal_{\tilde{\Sigma}}^{\pm}$ contains a disk $D^{\pm}$ in the sphere at 
infinity of radius bounded below by a constant that only depends on the quasiconformal dilatation and $U$ (the actual constant doesn't matter for the proof.)  

It also follows from \cite{Eps86} that the normal exponential map $\eta$ and Gauss maps $\Gcal_{\tilde{\Sigma}}^{\pm}$ define an embedding $\overline{\eta}_U$ from $U \times [-\infty,\infty]$ into the closed Poincare ball $\mathbb{H}^3 \cup \partial_\infty \mathbb{H}^3$. Note that $D^\pm$ is contained in the image of $\overline{\eta}_U$.  The convex hulls of $\overline{D}^{\pm}$ define solid hemisphere regions $H^{\pm}$ in $\HH$.  By making $D^\pm$ smaller if necessary, we can assume that the $H^{\pm}$ are disjoint from $U$, and thus the rest of $\tilde{\Sigma}$, since $\eta$ is a diffeomorphism.    

Note that $D^{\pm}$ are disjoint from the limit set of $\tilde{\Sigma}$, or equivalently that $D^{\pm}$ are contained in the domain of discontinuity for the action of $\pi_1(\Sigma)$ on $\partial_\infty \HH$.  This is true because the $\Gcal_{\tilde{\Sigma}}^{\pm}$ are $\pi_1(\Sigma)$-equivariant homeomorphisms and $U$ is disjoint from all of its translates under the covering action of $\pi_1(\Sigma)$ on $\mathbb{H}^3$. 

The regions $H^{\pm}$ thus both project to open regions outside of the convex core of $M$. Since they project to different ends of $M$, this shows that $M$ has no 
degenerate ends. Thurston \cite{Thu79} proved that each end of $M$ is either geometrically finite or 
has a neighborhood contained in the convex hull.  It follows that both ends of $M$ are geometrically finite, and in the case that $M$ contains no {\ap}s,  $M$ must also be quasi-Fuchsian. \qed 

\begin{rem} \label{rubinsteinexample} 
Rubinstein \cite{Rub05} constructed examples of essential immersed surfaces $\Sigma$ in the figure eight knot complement with principal curvatures less than or equal to 1 in magnitude, for which the associated Kleinian surface group has an accidental parabolic. The accidental parabolics correspond to embedded horocycles that lie on the surfaces $\Sigma$. Although the $\Sigma$ are not minimal, it seems possible that they are homotopic to minimal surfaces with $\lambda_0\le 1$. 
If this is the case, then there would be elements of $\Bcal_0$ for which the corresponding Kleinian groups had accidental parabolics.  
\end{rem} 
If there are {\ap}s in the corresponding {\ksg} of $M\in \Bcal_0$, we will take advantage of the explicit metric (see \eqref{gausscoordinates}) on $M$ to
prove that $M$ is geometrically finite. A key observation we will use is the following: the tangent vectors to $\Sigma$ that do not lie along principal directions 
with {\pc}s $\pm 1$ expand exponentially under the normal exponential map.

\noindent \textit{Proof of Theorem A1:} 

 The topology of the convex core of $M$ is determined by the following information \cite{BCM12}: there is a multicurve $C$ 
 on $S \times \{1\}$ each component of which corresponds to an upward cusp in $M$ and an {\ap} in the {\ksg}. Each component of the complement 
 of $C$ in  $S \times \{1\}$ either corresponds to an upward geometrically finite end or an upward degenerate end.  The downward end 
 is described in a similar way. The only compatibility condition is that the multicurve $C'$ for the downward end not contain any curves homotopic to curves in $C$.  The convex core of $M$ is then homeomorphic to $\Sigma \times [0,1]$ with $C \times \{1\}$ and $C' \times \{0\}$ removed, along with any connected components of the complement of $C \times \{1\}$ in $S \times \{1\}$ and $C' \times \{0\}$ in $S \times \{0\}$ corresponding to degenerate ends.  
 
We now give a more detailed description of $M$ in the case that $M \cong S \times \mathbb{R}$ for $S$ a closed surface, following \cite{Min10}. Let $Q$ denote the union of the $\epsilon$-Margulis tubes of rank-1 cusps of $M$, where $\epsilon$ is taken small enough that all of the Margulis tubes are disjoint (each such Margulis tube corresponds to an accidental parabolic). Then if $M_0= M - Q$, we can choose a compact three-dimensional submanifold $K$ with boundary of $M_0$ whose inclusion in $M_0$ is a homotopy equivalence.  Moreover $K$ can be chosen so that $\partial K$ meets the boundary of each rank-1 cusp of $Q$ in an essential annulus. The ends of $M_0$ are exactly connected components of $M_0 - K$, which are in one-to-one correspondence with connected components of $\partial K - \partial Q \cap K$. Let $E$ be an end of $M_0$, i.e. some connected component of $M_0 - K$. Then we can choose a homeomorphism $\Phi_E$ from $\Sigma_0' \times [0,\infty)$ onto the closure of $E$, for $\Sigma_0'$ some compact subsurface of $\Sigma$.   The surface $\Sigma_0'$ corresponds to a connected component of the complement of the $C$ or of the $C'$ from the previous paragraph.

  We know that any point  in $\partial E$ that is not contained in $\partial K$  must  be contained in the $\epsilon$-thin part of $M$, and so this is true of any point in the boundary of $\Phi_E(\Sigma_0' \times \{t\})$ for $t\geq 0$. For any boundary component of $\Sigma_0'$,  we can thus choose an essential loop $\gamma_t$ with length less than $\epsilon$ homotopic to and intersecting a loop in the boundary of $\Phi_E(\Sigma_0' \times \{t\})$. The homotopy class of $\gamma_t$ corresponds to an accidental parabolic.

  Recall that the normal exponential map 
$\eta:\Sigma \times \mathbb{R} \rightarrow M$ for $\Sigma$ is a diffeomorphism and gives global coordinates for $M$. In these coordinates,  Uhlenbeck \cite{Uhl83} showed that the metric can be written as: 
\be\label{gausscoordinates} 
e^{2v(x)}(\cosh(r) \text{I} + \sinh(r) e^{-2v(x)} A(x))^2 + dr^2,
\ene
\noindent where  $A(x)$ is the {\sff} of $\Sigma$ at $x$ and $\text{I}$ is the identity matrix.  Here $x$ is a local isothermal coordinate on $\Sigma$, so that the hyperbolic metric on $\Sigma$ at $x$ is given by $e^{2v(x)} \text{I}^2$.   

Note that the existence of global normal coordinates on $M$ relative to $\Sigma$ as in (\ref{gausscoordinates}) implies that for any sequence of points in an upward end $E$, the $t$-coordinate of the $\Phi_E$-preimages tends to infinity if and only if the signed distance $r$ to $\Sigma$ tends to infinity.  It is clear from (\ref{gausscoordinates}) that any tangent vector to $\Sigma$ gets exponentially expanded as $r$ tends to infinity unless it is tangent to a principal direction with {\pc} $-1$. Note that the function $f=|A|^2$ attains its maxima precisely at the points where the {\pc}s are $\pm 1$ and that $f$ is real analytic on $\Sigma$ 
since $A$ is the real part of a holomorphic quadratic differential (see e.g. [Proof of Lemma 4.1, \cite{WW20}]). Therefore the set of points $\Gamma$ where $\Sigma$ has a {\pc} equal to $\pm 1$ has connected components that are either isolated points or embedded graphs.

In the normal exponential coordinates $\eta$, since the length of $\gamma_t$ is bounded from above by $\epsilon$, we have that $\gamma_t \subset \eta (\Sigma \times [r_t,r_t+\epsilon])$ for $r_t$ tending to infinity as $t$ tends to infinity. By projecting $\gamma_t$ to the $\Sigma$-factor under $\eta$ we obtain a curve $\hat{\gamma}_t$ in $\Sigma \times \{0\}$ homotopic to $\gamma_t$. Let $z_t$ be a function from $\hat{\gamma}_t$ to $[r_t,r_t +\epsilon_t]$ defined by $\gamma_t(s) = \eta(\hat{\gamma}_t(s), z_t(s)$ for some variable $s$ parametrizing both $\gamma_t(s)$ and $\hat{\gamma}_t(s)$. By the exponential expansion of tangent vectors not tangent to a principal direction with principal curvature $-1$, we know that for every 
$\epsilon'>0$ the $\epsilon'$-neighborhood of $\Gamma$ will contain $\hat{\gamma}_t$ for sufficiently large $t$. Since the $\hat{\gamma}_t$ are essential, we know that for $t$ large enough they must be $\epsilon'$-close to the components of $\Gamma$ that are embedded graphs (as opposed to the components that are isolated points.) 

For $n$ large enough there is a well-defined projection of the $\hat{\gamma}_t$ to $\Gamma$. By projecting to $\Gamma$ we thus obtain curves $\overline{\gamma}_t$ in the same homotopy class as $\hat{\gamma}_t$. By eliminating backtracking in $\overline{\gamma}_t$, we obtain a concatenation $\mathcal{C}$ of edges in $\Gamma$. We claim that each of these edges travels everywhere tangent to a principal curvature -1 principal direction, provided $n$ was taken large enough.  To see this, suppose for contradiction that for infinitely many $n$ there were some edge $I$ of $\Gamma$  in $\mathcal{C}$, such that $I$ contained a point $p$ whose tangent vector to $I$ was not tangent to a principal curvature -1 principal direction.  Then this will also be the case for some sub-interval $I'$ of $I$ containing $p$.  

Let $R(\epsilon')$ be the rectangular region containing $I'$ defined by taking the image of $I' \times (-\epsilon',\epsilon')$ under the normal exponential map of $I'$ in $\Sigma$. Then for any $L$, provided $\epsilon'$ is chosen small enough and $t$ is chosen large enough, the following will hold: Let $c$ be a curve joining the two ``thin"  boundary components of $R(\epsilon')$-- i.e., the boundary components corresponding to $\partial(I') \times (-\epsilon',\epsilon')$, which by making $\epsilon'$ small we can take to have much smaller length than $I'$.  Then the image of $c$ under the map $\eta(c(\cdot),f(c(\cdot))$ has length at least $L$, for any curve $c$ as above and any smooth function $f> t$ defined on $c$.  

For $t$ large, we know that the curves $\hat{\gamma}_t$ intersected with $R(\epsilon')$ will contain curves $c$ as in the previous paragraph.  This implies that $\gamma_t(\cdot) = \eta(\hat{\gamma}_t(\cdot ), z_t(\cdot))$ will have length tending to infinity with $t$, which contradicts the fact that the $\gamma_t$ have uniformly bounded length.  This proves that the concatenation of edges $\mathcal{C}$ obtained from $\overline{\gamma}_t$ is a line of curvature.  The same argument also shows that $\overline{\gamma}_t$ is contained in a neighborhood of $\mathcal{C}$ of radius tending to zero as $t$ tends to infinity.  This is because, for a  vertex of $\Gamma$ where $\overline{\gamma}_t$ veers off of $\mathcal{C}$, the edge along which $\overline{\gamma}_t$ veers off will fail to be tangent to a principal curvature $-1$ direction in some neighborhood of $v$.  $\overline{\gamma}_t$ can therefore venture only $\delta(t)$ far along this edge before backtracking, where $\delta(t)$ tends to zero as $t$ tends to infinity.

It follows that for $t$ sufficiently large each boundary component $\gamma_t$ of $\Sigma_0' \times \{t\}$ will map to a small neighborhood of an embedded line of curvature $\overline{\gamma}_t'\subset \Gamma_0$ under normal projection to $\Sigma$ using as above the normal exponential coordinates given by $\eta$.  By passing to  subsequence of times $t_n \rightarrow \infty$, we can assume that $\overline{\gamma}_{t_n}'$ is independent of $t_n$ for each boundary component of $\Sigma_0'$.  The $\overline{\gamma}_{t_n}'$ bound a surface in $\Sigma$ isotopic to $\Sigma_0'$: choose some point $p$ in the interior of this surface that is not contained in $\Gamma$.  

The normal projection of $\Phi(\Sigma_0' \times \{t_n\})$ down to $\Sigma$ for $t_n$ large enough will then contain a fixed neighborhood $U$ of p on which the {\pc}s are strictly less than $1$.  Consequently for every $T>0$ there is a geodesic segment that begins normal to $U$, has length longer than $T$, and has an endpoint contained in $E$. We can then lift $U$ to the universal cover and apply the argument from the proof of Theorem A above to produce points on upward normal geodesic rays from $U$ that are not contained in the convex core of $M$, which as before proves that $E$ is non-degenerate.  It follows that $M_0$ has only non-degenerate ends, and that $M$ is geometrically finite.

\qed

An immediate consequence is that elements of $\Bcal_0$ which have no accidental parabolics are quasi-Fuchsian, which is Corollary ~\ref{cor1}.  

Although the statements of Theorems A and A1 are for closed surfaces $S$, essentially the same proof applies to prove them for punctured surfaces $S'$ of finite type.  This will be important for the construction of the compactification in the next section, and we write it as a corollary.  

 \bcor \label{cusp}
Let $M_1$ be a {\htm} diffeomorphic to $S'\times \R$, where $S'$ is a surface that carries a finite volume hyperbolic metric. If $M_1$ admits 
an incompressible minimal surface $\Sigma'$ (homotopic to $S' \times \{0\}$) with $|\lambda(\Sigma')|\le 1$ then $M_1$ is geometrically finite.  
\ecor

\begin{proof} 

First, assume that the only parabolics correspond to the cusps of $\Sigma'$.  Then we can lift $\Sigma'$ to a properly embedded disc in $\mathbb{H}^3$.  Note that there is a point on the minimal surface where the principal curvatures are strictly less than one in absolute value.  In the closed case this followed from the Riemann-Roch theorem, but here it also follows from the fact that the holomorphic quadratic differential whose real part is the second fundamental form decays rapidly at any cusp. 

To see this, lifting the holomorphic quadratic differential to the universal cover $\mathbb{H}^2$ of $\Sigma'$ in its hyperbolic metric we obtain a weight 4 modular form $\phi(z) dz^2$. Given a cusp in $\Sigma'$ we can assume that it corresponds to the Mobius transformation $z\mapsto z+1$ with fixed point $\infty \in \partial_{\infty} \mathbb{H}^2$.  That the principal curvatures are bounded by one implies that $y^2 \phi(z)$ is bounded above in absolute value as y tends to infinity.  Since $\phi(z)$ is holomorphic and invariant under $z \mapsto z+1$, it has a Fourier expansion $\sum_{n=1}^\infty a_n e^{2\pi i n z }$ \cite{s73}.
  \noindent This implies that $|\phi(z)|$ decays exponentially fast at the cusp. Note that this also implies that the set of points with principal curvatures $\pm 1$ is contained in a compact set, and is homeomorphic to a union of isolated points and embedded graphs as above.   

Thus we can run the same argument as the first part of the proof of Theorem A to produce neighborhoods of both connected components of the complement of $\Sigma'$ not contained in the convex core.  It then follows by Thurston that $\Sigma'$ is quasi-fuchsian.  

In the case that there are accidental parabolics, we still have decompositions of $\Sigma'$ into subsurfaces $\Sigma_0'$ as in the first two paragraphs of the proof of Theorem A1 above, one such decomposition for the upward and one for the downward end \cite{Min10}[pgs. 11-12]. (We thus argue first for the upward end, using the decomposition of $\Sigma'$ corresponding to the upward end, and then for the downward end). The only difference in that description are that the components $Q$ of the $\epsilon$-thin part may correspond to punctures of $\Sigma'$ in addition to accidental parabolics, and that the surfaces $\Sigma_0'$ may have cusps. But the same argument applies to produce arbitrarily long normal geodesic segments from these $\Sigma_0'$ with an endpoint not contained in the convex core, thus ruling out degenerate ends.    
\end{proof} 
There is also a notion of a closed surface in a hyperbolic 3-manifold being {\qf}. A closed surface $S$ of genus $g \ge 2$ in a complete {\htm} $N$ is called {\qf} if a lift of the inclusion of the 
universal covers is a quasi-isometry. This is equivalent to $S$ being $\pi_1$-injective and the cover of $N$ corresponding to $\pi_1(S)$ being quasi-Fuchsian. 
A result of Thurston (proved in  \cite{Lei06}) states that a closed surface of genus at least 2 in a complete {\htm} is {\qf} if its {\pc}s are strictly less than 1 in magnitude.  
The proof of Theorem ~\ref{main1} generalizes this to the following:
 \bcor\label{qf}
If $N$ is a complete {\htm}, and $S$ is a closed surface in $N$ such that the {\pc}s are less than or equal to one in magnitude and strictly 
less than one in magnitude at some point, and if the Kleinian group corresponding to $S$ contains no accidental parabolics, then $S$ is {\qf}. 
 \ecor 
 \bp
 We  claim that $S$ is $\pi_1$-injective, or equivalently that the lift $\tilde{S}$ of $S$ to $\HH$ is a disc. Then the {\pc}s $\lambda(\tilde{S})$ of $\tilde{S}$ satisfy that $|\lambda(\tilde{S})| \le 1$ and there 
exists some point $\tilde{p} \in \tilde{S}$ such that $|\lambda(\tilde{p})| < 1$. If $\tilde{S}$ were not homeomorphic to a disk, then taking a closed geodesic in $\tilde{S}$ in its induced metric and applying the argument in the proof of \cite{Eps84}[Theorem 3.4] would give a contradiction as follows: Epstein showed that the 
hyperbolic cosine of the distance from the starting point of a curve in $\HH$ with geodesic curvature less than or equal to $1$ in absolute value (such 
as a geodesic on $\tilde{S}$ in its induced metric) is convex along that curve, and that therefore such a curve cannot return to its starting point. 
It follows that $\tilde{S}$ must be a disc. The argument can then proceed as in the proof of Theorem A above to show that $S$ is {\qf}.
 \ep
 \section{Compactifying Weakly Almost Fuchsian Space} \label{compactificationsection} 
 
 In this section, we construct the compactification $\overline{\Bcal}_0$ of the space of unmarked weakly almost Fuchsian manifolds $\Bcal_0$. Our 
 compactification extends the Deligne-Mumford compactification of the moduli space of {\RS}s, and is analogous to the compactification defined 
 by Canary-Storm of the space of unmarked {\ksg}s \cite{CS12}.

 \noindent \textit{Proof of Theorem B:} 

We first describe the points we add to the space of weakly almost Fuchsian manifolds in order to form the compactification $\overline{\Bcal}_0$.  This will involve defining a potentially larger space $\overline{\Bcal}_0'$. Next we define what it means for a sequence to converge in $\overline{\Bcal}_0'$, and show that $\Bcal_0$ is sequentially precompact in $\overline{\Bcal}_0'$. We then construct  $\overline{\Bcal}_0$ by taking the closure of $\Bcal_0$ in $\overline{\Bcal}_0'$.


Our construction utilizes a triple of data attached to the unique closed {\ms} in a {\waf} manifold. Uhlenbeck \cite{Uhl83} showed that any triple $(g_\sigma, e^{2u},\alpha)$ of a hyperbolic metric on a closed surface, conformal factor, and {\hqd} on $(\Sigma,g_\sigma)$ 
 that satisfies the Gauss equation and has {\pc}s less than or equal to one gives a unique complete hyperbolic structure $M$ on $\Sigma \times \mathbb{R}$ that contains a {\ms} with {\sff} given by the real part $\Re(\alpha)$ of $\alpha$, {\pc}s no more one in magnitude, and induced metric $e^{2u} g_\sigma$. This is the  unique closed {\ms} in $M$. 
 
 Furthermore, under the {\pc} bound $\lambda_0 \le 1$, the solution $u$ for the Gauss equation is unique. That $\alpha$ 
 be holomorphic is equivalent to the {\sff} it defines satisfying the Codazzi equations, provided the surface is minimal. When we say that a triple 
 $(g_\sigma,e^{2u},\alpha)$ satisfies the Gauss equation and has principal curvatures less than or equal to one, we mean that a {\ms} in a {\htm} with induced 
 metric $e^{2u} g_\sigma$ and {\sff} given by $\Re(\alpha)$ has this property if it exists. We know a posteriori that a {\ms} with this data exists and is contained in a complete hyperbolic 3-manifold by \cite{Uhl83}.  By Corollary ~\ref{cusp}, this hyperbolic 3-manifold is geometrically finite.

For each genus $g>1$, we define $\overline{\Bcal}_0'=\overline{\Bcal}_0'(g)$ as follows.  A point in $\overline{\Bcal}_0'$ will be given by a finite set of (possibly punctured) Riemann surfaces $\Sigma_1,..,\Sigma_n$ with unmarked conformal structures we denote by $\sigma_1,..,\sigma_n$, smooth functions $e^{2u_1},..,e^{2u_n}$ (which we view as conformal factors) on the $\Sigma_i$, and holomorphic quadratic differentials $\alpha_1,..,\alpha_n$ on the respective Riemann surfaces $(\Sigma_1,\sigma_1),..,(\Sigma_n,\sigma_n)$.  A point in the compactification also contains the information of an identification of the cusps of the $\Sigma_i$, so that when we glue compact cores of the $\Sigma_i$ according to this identification we obtain a closed orientable surface of genus g.  Finally we require that the $(\sigma_i,e^{2u_i},\alpha_i)$ satisfy the Gauss equation and have principal curvatures less than or equal to one in magnitude.

Each point in $\overline{\Bcal}_0'$ thus gives a finite disjoint union of cusped weakly almost Fuchsian manifolds $M_1,..,M_n$, that are geometrically finite and contain essential minimal surfaces $\Sigma_i$ with principal curvatures smaller than or equal to one, as above. (Although Uhlenbeck's results are stated for $\Sigma$ a closed surface, the construction of a unique hyperbolic 3-manifold from the triple of a conformal structure, conformal factor, and holomorphic quadratic differential satisfying the Gauss equation and having principal curvatures less than or equal to $1$ works equally for $\Sigma$ a surface that carries a finite volume hyperbolic metric). In what follows we will pass freely between the conformal structures $\sigma_i$ and the hyperbolic metrics $g_{\sigma_i}$ that uniformize them, which are unique up to isometry.  Also, rather than referring to all of the $\Sigma_i$ separately, we will denote a point in $\overline{\Bcal}_0'$ by the triple $((\Sigma,\sigma), e^{2u}, \alpha)$, where $(\Sigma,\sigma)$ is a possibly disconnected unmarked Riemann surface, and $e^{2u}$ and $\alpha$ are respectively a function and holomorphic quadratic differential on $\Sigma$. In referring to a point in $\overline{\Bcal}_0'$ this way we take the pairings of the cusps to be implicitly given.  We comment that for any $((\Sigma,\sigma), e^{2u}, \alpha)$ in  $\overline{\Bcal}_0'$, the number of connected components of $\Sigma$ is finite and bounded by a constant depending on $g$.  This follows from the fact that the connected components of $\Sigma$ are all hyperbolic surfaces, and by gluing compact cores based on the pairings of the cusps we obtain a closed genus g surface.

We now define what it means for a sequence to converge in $\overline{\Bcal}_0'$. We say that a sequence $(\sigma_k, e^{2u_k},\alpha_k) \subset \overline{\Bcal}_0'$ converges to $(\sigma, e^{2u},\alpha)$ if the following holds. First we require that $\sigma_k$ converges to $\sigma$ in the Deligne-Mumford compactification of the moduli space of Riemann surfaces of genus $g$.  Next, let $\Sigma$ be the possibly disconnected Riemannian surface whose metric is given by  $e^{2u} g_{\sigma}$, and let $\Sigma_k$ be the Riemannian surfaces corresponding to $(g_{\sigma_k}, e^{2u_k},\alpha_k)$ in the same way. Then there exist possibly disconnected compact subsurfaces $C_k$ exhausting $\Sigma$ and smooth maps $\Phi_k:C_k \rightarrow \Sigma_k$ so that for all large enough $k$: 

\begin{itemize}
\item The $\Phi_k$ map any two essential loops in $\Sigma$ that are homotopic into respective cusps that are paired in $\Sigma$ to either homotopic curves in $\Sigma_k$ or loops homotopic into paired cusps in $\Sigma_k$.  
	\item $\Phi_k$ is a diffeomorphism onto its image.   
 \item As $k$ tends to infinity the maps $\Phi_k$ are smoothly converging to isometries, and the pullbacks of the $\alpha_k$ under $\Phi_k$ are smoothly converging to $\alpha$ on compact sets.
	\end{itemize}

We now check that $\Bcal_0$ is sequentially precompact in $\overline{\mathcal{B}}_0'$. Take a sequence of  hyperbolic 3-manifolds $M_k$ corresponding to elements of $\Bcal_0$.  Let $\Sigma_k$ be the unique closed embedded {\ms} in $M_k$, $\sigma_k$ be the {\cs} of its induced metric, and $\alpha_k$ be the {\hqd} in $(\Sigma_k, \sigma_k)$ 
that encodes the {\sff} of the {\mi}. We write the induced metric on $\Sigma_k$ as the {\hym} $g_{\sigma_k}$ multiplied by a conformal factor $e^{2u_k}$.

We can pass to a subsequence, which by abuse of notation we also denote by $(g_{\sigma_k},e^{2u_k},\alpha_k)$, so that the unmarked conformal structures $\sigma_k$ converge to a point in the Deligne-Mumford compactification of the moduli space. This point is given by a disjoint union of cusped surfaces $\overline{\Sigma}_1,..., \overline{\Sigma}_n$, which we denote by $\overline{\Sigma}$. The convergence of the conformal structures implies smooth convergence of the uniformizing hyperbolic structures on compact sets \cite{b74}.  There consequently exist smooth maps $\Phi_k:\overline{\Sigma} \rightarrow \Sigma_k $ that smoothly converge to isometries on compact sets.

 Next, we claim that we can pass to a subsequence of the $u_k$ whose pullbacks by the $\Phi_k$ converge smoothly on compact sets.  We first show that the $u_k$ are  uniformly bounded in $L^\infty$. That  $u_k \le 0$ follows from the {\maxp} as in \cite{Uhl83}. Furthermore the {\pc} condition $\lambda_0 \le 1$ implies that the Gaussian curvatures are $-1-\lambda_0^2 \ge -2$. By the conformal change equation the Gaussian curvature is given by $e^{-2u_k}(-1-\Delta_{\sigma_k} u_k)$, where $\Delta_{\sigma_k}$ is the Laplace operator for the {\hym} $g_{\sigma_k}$. Therefore we have 
 $e^{-2u_k}(-1-\Delta_{\sigma_k} u_k) \ge -2$, and we deduce by the {\maxp} that $u_k \ge \frac{-\ln(2)}{2}$.  Because the $u_k$ are bounded in $L^\infty$, satisfy a second order elliptic equation, and the $\Phi_k$ converge smoothly to isometries on compact sets,  by standard elliptic theory we can pass to a subsequence of the $u_k$ for which the pullbacks under $\Phi_k$ converge smoothly on compact sets.

Finally the $\alpha_k$'s have uniformly bounded $L^\infty$-norm, measured in the hyperbolic metrics $g_{\sigma_k}$. This is because the principal curvatures are bounded above by 1 in absolute value and at each point the principal curvatures $\pm \lambda$ are equal to $\pm |\alpha_i|e^{-2u_i}$, where the norm of $\alpha_i$ is measured in the hyperbolic metric $g_{\sigma_k}$ (see \cite{Uhl83}[pg. 157]), and the conformal factors $u_k$ satisfy $-\frac{ln(2)}{2} \leq u_k \leq 0$.  Therefore we can pass to a further subsequential limit so that the pullbacks of the $\alpha_k$ by $\Phi_k$ are smoothly converging to a holomorphic quadratic differential $\alpha$ on compact subsets.  We thus obtain a triple in $\overline{\Bcal}_0'$ to which a subsequence of our original sequence converges in the sense just defined.  Two cusps are paired if they correspond to the same curve that gets pinched off in the convergence of the $\sigma_k$ to $\sigma$ in the Deligne-Mumford compactification, or else if they correspond to cusps of $\sigma_k$ that were already paired for infinitely many $k$ in our subsequence.

 We define $\overline{\Bcal}_0$ to be the closure of $\mathcal{B}_0$ in $\overline{\mathcal{B}}_0'$: i.e., the set of all points in $\overline{\mathcal{B}}_0'$ that are limits of sequences in $\mathcal{B}_0$. As we have just shown, $\Bcal_0$ is sequentially precompact in $\overline{\mathcal{B}}_0'$, so $\overline{\mathcal{B}}_0$ is sequentially compact (one can check this directly by a diagonal argument.) 

  As we verify in the proof of the next proposition, each Kleinian group corresponding to a connected component of the hyperbolic 3-manifold for an element of  $\overline{\Bcal}_0$ arises as an algebraic limit of remarked representations for the $M_k$, possibly restricted to subsurfaces.  One could thus likely topologize $\overline{\Bcal}_0$ by viewing it as a space of Kleinian groups equipped with maps of minimal surfaces into their quotients.  This would be closer to the  Canary-Storm approach to compactifying spaces of Kleinian groups (\cite{CS12}).  

  To complete the proof of Theorem ~\ref{main2}, we need to analyze more closely the convergence of the $M_k$ to the disjoint union $\sqcup_{i=1}^m \overline{M}_i$.  For the sake of simplicity and since it is all that is needed for the applications, we restrict in the following proposition to weakly almost Fuchsian points of $\overline{\Bcal}_0$.  


\bpo \label{boundarycontinuity} 
Let $\{M_k\} \subset \Bcal_0$ be a sequence of weakly almost Fuchsian manifolds of topological type $\Sigma \times \mathbb{R}$, and suppose that the $M_k$ converge to a point $\sqcup_{i=1}^m \overline{M}_i$  in the compactification.  Then there is a uniform bound (depending on  $\sqcup_{i=1}^m \overline{M}_i$)  on the volume of the convex cores of the $M_k$.

\epo
 \bp	  

Let $c_1,..,c_\ell$ be the simple closed curves on the {\ms} $\Sigma_k \subset M_k$ which become nodes in the limit. We take these curves to be the unique geodesics in their homotopy class on $\Sigma_k$ in its induced metric.  We claim that normal neighborhoods in $M_k$ of each connected component $C_i(k)$ of the complement of the disjoint  union of the $c_j$ are converging to the $\overline{M}_i$ on compact subsets, $i=1,..,m$ (we implicitly choose some consistent marking and identification of all of the $\Sigma_k$ so that this makes sense). By normal neighborhood we mean the image of some subset of the form $C_i(k) \times (-L,L)$ of the normal bundle to $C_i(k)$ under the normal exponential map.

 More precisely, for each $\overline{M}_i$ there exists a map $h_k^i: \overline{M}_i\rightarrow M_k$ which is a homotopy equivalence onto its image, whose image is a normal neighborhood of the complementary region $C_i(k)$, and which restricted to any compact three-dimensional submanifold with boundary of $\overline{M}_i$ is a diffeomorphism onto its image for large enough $k$.  These maps are defined as follows. For each $\epsilon>0$ smaller than the Margulis constant and the shortest closed geodesic on $\overline{M}_i$, there is a retraction $r_{\epsilon,k}$ from $\overline{M}_i$ to its $\epsilon$-thick part (see the discussion on the truncated convex core in \cite{McM99}[Section 4].)  Here $\epsilon=\epsilon(k)$ will tend to zero as $k$ tends to infinity, in such a way that the following holds: there is a diffeomorphism $\Phi_{\epsilon,k}$ between a surface with boundary that contains the intersection of $\overline{\Sigma}_i$ with the $\epsilon$-thick part of $\overline{M}_i$,  and a compact subsurface of $\Sigma_k$ obtained by taking a complementary region to a set of curves each homotopic to one the $c_j$.  Moreover, for any compact submanifold $K$ of $\overline{\Sigma}_i$ and for $k$ large enough, $\Phi_{\epsilon,k}$ restricted to $K$ is as $C^\infty$ close as desired to being an isometry onto its image.  To define $h_k^i$ we apply $r_{\epsilon,k}$, take the preimage under the normal exponential map for $\overline{\Sigma}_i$ to obtain a point in $\overline{\Sigma}_i \times \mathbb{R}$, apply $\Phi_{\epsilon,k}$ to the first factor $\overline{\Sigma}_i$ while keeping the second factor $\mathbb{R}$ fixed, and then take the image under the normal exponential map for $\Sigma_k$.  
 
 Note that the maps $h_k^i$ are smoothly converging to isometries as $k$ tends to infinity. This follows from the explicit formula 
 \eqref{gausscoordinates} for the metric on the normal neighborhood of the minimal surfaces $\Sigma_k$ with $\lambda_0(\Sigma_k)\le 1$, and the fact that the metrics on the $C_i(k)$ are smoothly converging to $e^{2\overline{u}_i}$ times the hyperbolic metric on $\overline{\Sigma}_i$, and similarly for the associated holomorphic quadratic differentials. 
 
 Fix some $i$. It then follows, for a choice of 
 baseframes $q_k=q_k(i)$ for the $M_k$, each the $h_k^i$-image of some fixed frame $q$ over a point in $\overline{M}_i$, that $(M_k,q_k)$ converges geometrically to $(\overline{M}_i,q)$. 
 Here the choice of baseframes amounts to, in the limit, throwing out the complement of $\pi_1(C_i(k))$ in $\pi_1(M_k)$.  That we have geometric convergence can be seen as follows. To check geometric convergence, we need to check that $(\overline{M}_i,q)$ can be exhausted by compact sets $F_k$  that are mapped to $M_k$ by $C^1$ $1+\delta(k)$-quasi-isometries for $\delta(k)$ tending to zero as $k$ tends to infinity.  Here the $F_k$ are projections to $M_k$ of balls in the universal cover $\mathbb{H}^3$ centered at a lift of $q$ to $\mathbb{H}^3$, and we require that $q$ maps to $q_k$.  By the previous paragraph, a sequence of maps and $F_k$ satisfying these conditions can be constructed by restricting the maps $h_k^i$ to a sequence of compact sets exhausting $(\overline{M}_i,q)$.  
 
 We point out that to say that the $h_k^i$ are $1 +\delta(k)$ quasi-isometries, we are using the fact  that the normal exponential maps from the minimal surfaces $\Sigma_k$ and $\overline{\Sigma}_i$ are \textit{global} diffeomorphisms.  We also point out that the failure of this last property for hyperbolic 3-manifolds of type $\Sigma \times \R$ containing minimal surface with principal curvatures tending to but strictly greater than one is what prevents us from obtaining geometric convergence in the context of Section 5.

We moreover have that for all $i$, the Kleinian group $\Gamma_{C_i(k)}$ obtained by restricting the Kleinian group for $M_k$ to $C_i(k)$ converges algebraically, and therefore strongly, as $k\to \infty$ to a Kleinian group $\overline{\Gamma}_i$ such that $\overline{M}_i=\mathbb{H}^3/\overline{\Gamma}_i$.  Algebraic convergence follows from the fact that the maps  $h_k^i: \overline{M}_i\rightarrow M_k$ lift to maps to $\mathbb{H}^3 / \Gamma_{C_i(k)}$ that are homotopy equivalences and whose restrictions to any compact set are  $C^\infty$ converging to a local isometry. (Recall the definition of algebraic convergence in Section \ref{prelimksg}.)

Verifying the next claim is the key step in finishing the proof.  
 
\begin{claim} \label{keyclaim}
Fix a compact set $K\subset \overline{M}_i$.  Then the intersections with $K$ of the preimages under the $h_k^i$ of the convex cores $C_k$ of the $M_k$ are contained in the $C'$-neighborhood of the convex core of $\overline{M}_i$ for large enough $k$, for $C'$ any constant greater than the constant $C_0$ from Lemma \ref{density}.  
\end{claim}

 To prove the claim, we will use the fact that every point in the convex core of $M_k$ is near some closed geodesic (Lemma \ref{density}) together with the McMullen-Taylor curve straightening argument \cite{McM99} \cite{Tay97}.

 Let $\gamma$ be a nontrivial homotopy class of loop in $\Sigma_k$.  Then either (1) $\gamma$ is homotopic to a curve contained in some $C_i(k)$, or (2) we can decompose $\gamma$ minimally as a composition of homotopy classes of segments $\gamma_0,..,\gamma_{\ell-1}$ relative to their startpoints and endpoints, where each $\gamma_j$ has its interior contained in some $C_{i(j)}(k)$, and its startpoint and endpoint contained in respectively $c_{i^-(j)}$ and $c_{i^+(j)}$. We first handle the more difficult case of $\gamma$ in the second category. By a slight abuse of notation, for each $k$ we take $c_j=c_j(k)$  to be the unique geodesic on $\Sigma_k$ in the homotopy class of $c_j$, so that in particular the length of $c_{j}$ tends to zero as $k\to \infty$ and for any $\epsilon'>0$ the endpoints of the $\gamma_j$ are contained in the $\epsilon'$-thin part for $k$ sufficiently large.  We note that there is some ambiguity in the choice of a minimal decomposition of $\gamma$ into $\gamma_j$-- for example, we could modify it by composing the endpoint of $\gamma_j$ with $c_{i^+(j)}$ and the startpoint of $\gamma_{j+1}$ with  $c_{i^+(j)}=c_{i^-(j+1)}$ traversed in the other direction-- but that this will not matter for the arguments that follow. 


Following McMullen-Taylor, we now build a representative for $\gamma$ that is close to a geodesic, by modifying each of the $\gamma_j$ in turn and then combining them.  This construction will depend on choices of $m(k)$ and $\epsilon(k)$ so that $m(k)$ tends to infinity as $k$ tends to infinity, and so that not only $\epsilon(k)$ but $m(k)\epsilon(k)$ tends to zero as $k$ tends to infinity. The parameters $\epsilon(k)$ and $m(k)$ will also satisfy a couple of other conditions that we will specify below.

Fix some $\gamma_j$, and choose a segment $\gamma_j'$ in $\overline{\Sigma}_i \subset \overline{M}_i$, for $i=i(j)$, with endpoints contained in the boundary of the $m(k)\epsilon(k)/2$-thin part of the convex core of $\overline{M}_i$, so that $h_k^i(\gamma_j')$ is homotopic to $\gamma_j$ through segments on $\Sigma_k$ whose endpoints remain in the $m(k)\epsilon(k)$-thin part of $M_k$.   Since both endpoints of $\gamma_j'$ are contained in the convex core of $\overline{M}_i$, the unique geodesic in the relative homotopy class of $\gamma_j'$ is contained in the convex core of $\overline{M}_i$.  Denote by $\overline{\gamma}_j$ the image of this geodesic under $h_k^i$.   We can take $\epsilon(k)$ to be tending to zero in such a way that $h_k^i$ is $1+\delta(k)$ $C^2$-close to being an isometry when restricted to the $\epsilon(k)$-thick part of the convex core of  $\overline{M}_i$, where $\delta(k)$ tends to zero as $\epsilon(k)$ tends to zero. (Here we are using the fact that because $\overline{M}_i$ is geometrically finite, the $\epsilon'$-thick part of the convex core is compact for any $\epsilon'>0$.)   For $k$ large enough, we can thus take the geodesic curvature of  $\overline{\gamma}_j$ to be pointwise as small as desired away from the $\epsilon(k)$-thin part.

We modify $\overline{\gamma}_j$ as follows on the $m(k)\epsilon(k)$-thin region.  Let $J$ be a subsegment of  $\overline{\gamma}_j$ so that: 

\begin{enumerate} \label{twoconditions} 
	
	\item Both endpoints of $J$ are in the boundary of the $m(k)\epsilon(k)$-thin part.  

 \item The interior of J is contained in the $m(k)\epsilon(k)$-thin part.  
	
	\item $J$ has nontrivial intersection with the $\epsilon(k)$-thin part.  
	
\end{enumerate} 

We homotope $J$ to the unique geodesic segment $\overline{J}$ in its relative homotopy class.  Provided $k$ was chosen large enough, $\overline{J}$ also has non-zero intersection with the $\delta'(k)$ neighborhood of the $\epsilon(k)$-thin part for $\delta'(k)$ tending to zero as $k$ tends to infinity. To see this, if we assume for contradiction that we cannot find such a sequence of $\delta'(k)$, then for some uniform $\delta>0$ the geodesic segment $\overline{J}$ will stay at a distance of at least $\delta$ from the $\epsilon(k)$-thin part.  Since $\overline{J}$ is a geodesic, its preimage under $h_k^i$ has geodesic curvature tending to zero as $k$ tends to infinity.   It would then be possible to find a geodesic segment with the same endpoints in a neighborhood of $(h_k^i)^{-1}(\overline{J})$ of radius tending to zero as $k$ tends to infinity.  But by the fact that $J$ ventures into the $\epsilon(k)$-thin part there is a geodesic segment homotopic to $(h_k^i)^{-1}(\overline{J})$ relative to its endpoints that for large $k$ ventures into the $\delta$-neighborhood of the $\epsilon(k)$-thin part of $\overline{M}_i$.  These two geodesic segments are therefore distinct, which contradicts the fact that geodesic segments are unique in their relative homotopy class in negative curvature.  

 We have thus shown that $\overline{J}$ has non-zero intersection with the $\delta'(k)$-neighborhood of the $\epsilon(k)$-thin part, for $\delta(k)$ tending to zero as $k$ tends to infinity.  It follows that $\overline{J}$ meets $\overline{\gamma}_j$ at an angle  at each of its two endpoints that can be made as close to $\pi$ as desired by taking $k$ large. Here we are using the fact that $m(k)$ tends to infinity as $k$ tends to infinity, and that the geometry of any connected component of the thin part of $M_k$ approaches that of a standard cuspidal region as $k \to \infty$ (see (\cite{McM99}[pg. 14].)  We modify $\overline{\gamma}_j$ by replacing each subsegment $J$ satisfying (i),(ii), and (iii) above with $\overline{J}$ to obtain a new curve $\overline{\gamma}_j'$.  
 
 We now truncate the $\overline{\gamma}_j'$ by trimming off their ends in the following way: we remove the connected components of the intersection of $\overline{\gamma}_j'$ with the $m(k)\epsilon(k)$-thin part that contain the startpoint and endpoint of $\overline{\gamma}_j'$.  By abuse of notation we continue to denote this piecewise geodesic segment, which is a subsegment of the original $\overline{\gamma}_j'$, by $\overline{\gamma}_j'$.

 It remains to describe how to combine the $\overline{\gamma}_j'$ to obtain a nearly geodesic representative for $\gamma$. There are unique choices of relative homotopy classes of curves $\phi_j$ joining the endpoint of $\overline{\gamma}_j'$ to the startpoint of $\overline{\gamma}_{j+1}'$ so that when we concatenate the 
 $\overline{\gamma}_j'\phi_j$ we obtain a loop homotopic in $M_k$ to $\gamma$. 
 
 We replace the $\phi_j$ with the unique geodesic segments $\overline{\phi}_j$ in their relative homotopy classes and concatenate the $\overline{\gamma}_j'\overline{\phi}_j$ to obtain a loop $L$ homotopic to $\gamma$.  By the minimality of the decomposition of $\gamma$ into $\gamma_j$, we know that $\gamma_j$ approaches its endpoint from a different side of $c_{i^+(j)}=c_{(i+1)^-(j)}$ from which $\gamma_{j+1}$ approaches its startpoint (otherwise we could combine $\gamma_j$ and $\gamma_{j+1}$ into one subsegment, violating minimality). It therefore follows from the existence of global normal exponential coordinates from the $\Sigma_k$ that the length of $\overline{\phi}_j$ tends to infinity as $k$ tends to infinity.  Since the endpoints of $\overline{\phi}_j$ are both contained in the same connected component of the $m(k)\epsilon(k)$ thin part of the unique geodesic homotopic to $c_{i^+(j)}=c_{(i+1)^-(j)}$, the same is true for the interior of $\overline{\phi}_j$ for $k$ large enough that $m(k)\epsilon(k)$ is less than the Margulis constant.  (Here we are also using that $\overline{\phi}_j$ is homotopic relative to its endpoints to a curve contained in that component of the $m(k)\epsilon(k)$ thin part.)   Therefore $\overline{\phi}_j$ meets $\overline{\gamma}_{j}'$ and $\overline{\gamma}_{j+1}'$ at its endpoints at angles that tend to $\pi$ as $k$ tends to infinity. Here we are again using the fact that the geometry of any connected component of the thin part of $M_k$ approaches that of a standard cuspidal region as $k \to \infty$.
 
 Since $L$ can be divided into subsegments with lengths tending to infinity and geodesic curvature tending to zero as $k$ tends to infinity, and the angles at which these subsegments meet tend to $\pi$ as $k$ tends to infinity, the distance between $L$ and the unique geodesic in its homotopy class tends to zero as $k$ tends to infinity.  By how $L$ was constructed, for each $i$ the intersection of the preimage under $h_k^i$ of $L$ with the $m(k)\epsilon(k)$-thick part of $M_k$ is contained in the $\delta(k)$-neighborhood  of the convex core of $\overline{M}_i$, where $\delta(k)$ is independent of $\gamma$ and tends to zero as $k$ tends to infinity. In the case that $\gamma$ is homotopic to a curve contained in a single $C_i(k)$, we construct $L$ by taking the $h_k^i$-image of the corresponding geodesic in $\overline{M}_i$ and modifying the subsegments $J$ as above.  In the same way as before, the $L$ in this case is $\delta(k)$-close to the unique geodesic in its homotopy class and contained in the $\delta(k)$-neighborhood  of the convex core of $\overline{M}_i$, for $\delta(k)$ tending to zero as $k$ tends to infinity.  This together with the fact (Lemma \ref{density}) that every point in the convex core of $M_k$ is at a uniformly bounded distance from some closed geodesic gives the statement of Claim \ref{keyclaim}.

 It also implies that there are not pieces of the $m(k)\epsilon(k)$-thick part of the convex core of the $M_k$ that go off to infinity and result in the convex core losing volume in the limit.  To say this more precisely, for $k$ large enough every point in the $m(k)\epsilon(k)$-thick part of the convex core of $M_k$ is for some $i$ contained in the $C'$-neighborhood of the $h_k^i$-image of the $m(k)\epsilon(k)$-thick part of the convex core of  $M_k$, for $C'$ any fixed constant greater than the constant from Lemma \ref{density}.  This is because by the construction of the loops $L$, all points on $L$ that are a definite distance away from the $h_k^i$-images of the $m(k)\epsilon(k)$-thick parts of the $\overline{M}_i$ are contained in the $m(k)\epsilon(k)$-thin part of $M_k$.  As before, this is a consequence of the fact that geodesics joining any two points in a component of the $m(k)\epsilon(k)$-thin part of a hyperbolic 3-manifold, and that are homotopic relative to their endpoints to a curve in that component of the $m(k)\epsilon(k)$-thin part, stay in the $m(k)\epsilon(k)$-thin part their whole length provided $m(k)\epsilon(k)$ is smaller than the Margulis constant (Recall that the intersection of $L$ with the $m(k)\epsilon(k)$ thin part is a union of geodesic segments). 


For $k$ large enough that $m(k)\epsilon(k)$ is smaller than the Margulis constant, there are a uniformly bounded number of connected components of the $m(k)\epsilon(k)$-thin part of $M_k$, since every such connected component corresponds either to an accidental parabolic or one of the curves $c_j$ that get pinched off in the convergence of the $\Sigma_k$ to $\Sigma$ in the Deligne-Mumford compactification.  Each connected component of the $\epsilon'$-thin part of the convex core of $M_k$ has volume bounded above by some constant that tends to zero as $\epsilon'\to 0$.   Therefore for $k$ large the volumes of the convex cores of the $M_k$ are bounded above by the volume of the $C'$-neighborhood of the convex core of  $\overline{M}_i$, which is finite.

\ep 
We remark that the $c_j$'s in the proof are analogous to the shattering set considered in \cite{CS12}. 

Proposition ~\ref{boundarycontinuity} 
implies the Corollary ~\ref{cor3} in the introduction, which we restate here.   

\begin{cor} 
There exist $L$ and $\epsilon>0$ such that the volume of the {\cvc} and the {\Hd} of the limit set of any element of the set $\Bcal_0$ of weakly almost Fuchsian manifolds are bounded from above by respectively $L$ and $2-\epsilon$.  
\end{cor} 

\bp

Suppose for contradiction there is no uniform upper bound on the volume of the convex hulls of elements of $\Bcal_0$, and take a sequence of weakly almost Fuchsian $M_k$ whose convex cores have volumes tending to infinity. Then, passing to a subsequence, we can assume that the $M_k$ converge to some element of $\overline{\Bcal}_0$.  Proposition \ref{boundarycontinuity} then implies that the volumes of the convex cores of the $M_k$ have a uniform upper bound, which is a contradiction.

It follows from [Corollary A of \cite{BC94}] that for any $L$, there exists $\epsilon>0$ such that if the volume of the unit neighborhood of the {\cvc} of a 
geometrically finite infinite volume hyperbolic 3-manifold $M$ is bounded above by $L$, then the {\Hd} of the {\lmt} of $M$ is bounded above by 
$2-\epsilon$. The same argument as in the proof of Proposition ~\ref{boundarycontinuity} together with the argument of the previous paragraph shows that there is a uniform upper bound on the unit neighborhood of the convex core of any element of $\Bcal_0$, which finishes the proof.


\ep

We have thus proved Theorem ~\ref{main2} and Corollary ~\ref{cor3}, which answers Question ~\ref{qua1} and Question ~\ref{qua2} from the introduction. They also give criteria for a 
{\qfm} failing to be {\af}.
\qed

 \section{Beyond weakly almost Fuchsian Space} \label{finalsection}

We now give some applications and related results on {\ksg}s outside of the {\waf} case. As before we let $S$ be a closed surface and 
let $M$ be a complete {\htm} diffeomorphic to $S \times \mathbb{R}$. We first prove Theorem \ref{main4}. Before doing so, we recall the setup.   Let $M_n$ be a sequence of doubly degenerate hyperbolic 3-manifolds homeomorphic to $S \times \mathbb{R}$. Let $\Sigma_n$ be stable minimal surfaces in $M_n$ isotopic to $S \times \{0\}$, and suppose that the maximum principal curvatures of the $\Sigma_n$ tend to one.

 Then we claim that for some choice of markings $S \rightarrow M_n$ a subsequence $M_{n_k}$ of the $M_n$ converges algebraically on the complement of a multicurve $C$ in $S$.   To see this, note that as in the previous section we can pass to a convergent subsequence of the $\Sigma_n$ and their associated {\hqd}s and conformal factors to get a disjoint union $\sqcup_{i=1}^m (\overline{\Sigma}_i,\overline{\alpha}_i)$.  Since the {\pc}s are no more than $1$ in magnitude, we can construct hyperbolic structures $\overline{M}_i$ on $\overline{\Sigma}_i \times \mathbb{R}$ for which the $\overline{\Sigma}_i$ are  {\ms}s. We can then define maps $h_k^i$ as in the previous section, to show that on the complement of a multicurve a subsequence of the $M_k$ are algebraically converging to the disjoint union of the $\overline{M}_i$. What prevents us from obtaining geometric convergence in addition to algebraic convergence like in the previous section is that when the maximum principal curvature is larger than 1 the normal exponential map is no longer a global diffeomorphism.   

The following theorem implies Theorem \ref{main4} from the introduction, since assuming a lower bound on the injectivity radius algebraic and geometric limits are known to agree (\cite{McM99, Tay97}).  

\begin{thx} \label{alggeo} 
	
	For every connected component $S_0$ of the complement of $C$ in $S$, the restrictions to $S_0$ of the Kleinian groups $\Gamma(n_k)$ corresponding to the $M_{n_k}$ have the following property: every subsequential geometric limit of the $\Gamma(n_k)$ restricted to $S_0$ differs from the algebraic limit of $\Gamma(n_k)$ restricted to $S_0$.  
	
	\end{thx} 

\begin{proof} 
Note that each of the $\overline{M}_i$ is geometrically finite by Corollary \ref{cusp}, and that the homotopy equivalences $h_k^i$ are locally $C^\infty$ converging to isometries.  Suppose that for some complementary component $S_0$ of the multicurve $C$ corresponding to $\overline{M}_i$ it was the case that, up to a subsequence, the geometric limit was equal to $\overline{M}_i$ (that is to say, the algebraic and geometric limits agreed). Since $\overline{M}_i$ is geometrically finite, this would imply that $M_k$ contained points of arbitrarily large injectivity radius as $k$ tended to infinity. But this is impossible because there is a uniform upper bound on the injectivity radius of a doubly degenerate {\ksg} depending only on genus (see for instance \cite{Can96}).  Consequently for every complementary component $S_0$ every possible subsequential geometric limit must disagree with the algebraic limit, which completes the proof of Theorem \ref{alggeo}. 
\end{proof} 

\br
When $M$ is singly or doubly degenerate it is expected that $M$ contains a large number of closed {\ms}s. The existence of closed {\ms}s for some such $M$ is obtained in \cite{Cos21}, but the existence problem is still open in general.  
 \er  

We now prove 
Theorem ~\ref{main3}, restated below. We note that the proof of this theorem does not depend on the other theorems proved in this paper.
\bt
There exist {\qf} manifolds $M$ which contain unique closed {\ms}s $\Sigma$ each with {\pc}s strictly greater than 1 in absolute value at some point.  
\et
\bp

Take a path $\{M_t\}$ joining a Fuchsian manifold to a {\qfm} with multiple stable incompressible {\ms}s. Such examples were constructed in for instance \cite{HW15}. Let $t'$ be the greatest $t$ such that $M_t$ contains an incompressible {\ms} $\Sigma_{t}$ with {\pc}s less than or equal to 1. By \cite{HLT21}, 
we know that $\Sigma_{t'}$ is strictly stable, namely the bottom eigenvalue of the second variation operator 
\[
L= -\Delta_{\Sigma_{t'}} - |A|^2  +2 
\]
 \noindent of $\Sigma_{t'}$ is positive. The argument in appendix A of \cite{cg18}, which we reproduce in abridged form here, then shows that a neighborhood 
 of $\Sigma_{t'}$ has a mean-convex foliation. Let $\phi \in C^\infty(\Sigma_{t'})$ be an eigenfunction with bottom eigenvalue, which we can take 
 to be strictly positive, and let $N$ be the unit normal vector field to $\Sigma_{t'}$. Then if $F(x,t)$ is a variation of $\Sigma_{t'}$ with $F_t(x,0) = \phi \cdot N$ 
 and $\Sigma_{t'}(\tau)= F(\Sigma_{t'},\tau)$, then
\begin{equation} \label{meanconv} 
\frac{d}{d\tau}H_{\Sigma_{t'}(\tau)}|_{\tau=0} = L\phi= \lambda \phi >0.   
\end{equation}
 The $\Sigma_{t'}(\tau)$ for $\tau$ in some small interval about $0$ therefore give a mean-convex foliation (with respect to the outward normal vector) of a neighborhood of $\Sigma_{t'}$ in $M_t$.

 We now claim that for small enough $\epsilon$, $M_t$ has a unique closed {\ms} for $t \in [t',t'+\epsilon)$. For contradiction suppose not, and that there is a sequence of $t_n \searrow t'$  such that each $M_{t_n}$ has multiple closed {\ms}s. For $n$ greater than some large $N$, the implicit function theorem implies 
 that we can choose {\ms}s $\Sigma_{t_n}$ in $M_{t_n}$ converging to $\Sigma_{t'}$. Note that the $\delta$-neighborhood of each of the $\Sigma_{t_n}$ has a mean-convex foliation for $n>N$ and $\delta>0$ independent of $n$. If not, we could find a sequence of $\delta_n$ tending to zero and so that $\frac{d}{d\tau}H_{\Sigma_{s_n}(\tau)}(x_n)|_{\tau=\delta_n}=0 $ for $s_n \searrow t'$ as $n$ tends to infinity for some sequence of points $x_n$ on $\Sigma_{s_n}$.  By passing to a convergent subsequence of the $x_n$, we could find an $x$ on $\Sigma_{t'}$ so that $\frac{d}{d\tau}H_{\Sigma_{t'}(\tau)}(x)|_{\tau=0}=0$, which contradicts (\ref{meanconv}).

 Consequently there is some $\delta>0$ so that any other closed {\ms} $S_{t_n}$  in $M_{t_n}$ must contain a point at a distance of at least $\delta$ from $\Sigma_{t_n}$.  Fix a sequence of such $S_{t_n}$, and let $D_n$ be the maximum distance of a point on $S_{t_n}$ to $\Sigma_{t_n}$. We can pass to a subsequence of the $D_n$ that converges to some $D \in [\delta,\infty)$, by the fact that the convex cores of the $M_t$ are converging to the convex core of $M_{t'}$, and the fact that the ends of any {\qf} manifold have mean-convex foliations that serve as barriers (\cite{mp11}).  Note that the set of points $K_{t_n}$ at distance less than or equal to $D_n+\frac1n$ from $\Sigma_{t_n}$ is a compact submanifold for large enough $n$, that converges smoothly to the compact submanifold $K_{t'}$ with boundary of $M_{t'}$ consisting of points at distance less than or equal to $D$ from $\Sigma_{t'}$.   
 
 To summarize, the $S_{t_n}$ are contained in the $K_{t_n}$ and have points at distance from the boundary of $K_{t_n}$ tending to zero as $n$ tends to infinity.  The $K_{t_n}$ in turn are smoothly converging to $K_t$, which has strictly mean convex boundary  by the fact that $M_{t'} \in \Bcal_0$. But \cite[Theorem 3]{Whi10} says exactly that a sequence of minimal surfaces and convergent compact submanifolds as above is impossible, which gives a contradiction.

\ep
  

\bibliographystyle{amsalpha}
\bibliography{ref-gauss}
\end{document}